\documentclass{amsart}
\usepackage{amssymb,amsmath}

\newtheorem{theorem}{Theorem}
\newtheorem{lemma}{Lemma}
\newtheorem{example}{Example}
\newtheorem{problem}{Problem} 
\newcommand{\bt}{\begin{theorem}}
\newcommand{\et}{\end{theorem}}
\newcommand{\bl}{\begin{lemma}}
\newcommand{\el}{\end{lemma}}
\newcommand{\bex}{\begin{example}}
\newcommand{\eex}{\end{example}}
\newcommand{\bp}{\begin{problem}}
\newcommand{\ep}{\end{problem}}
\newcommand{\beal}{\begin{align*}}
\newcommand{\enal}{\end{align*}}
\newcommand{\beq}{\begin{equation}}
\newcommand{\eeq}{\end{equation}}
\newcommand{\benum}{\begin{enumerate}}
\newcommand{\eenum}{\end{enumerate}}
\newcommand{\ba}{\begin{array}}
\newcommand{\ea}{\end{array}}
\newcommand{\Z}{\ensuremath{\mathbf Z}}
\newcommand{\N}{\ensuremath{\mathbf N}}

\newcommand{\FF}{\ensuremath{\mathcal F}}
\newcommand{\GG}{\ensuremath{\mathcal G}}
\newcommand{\HH}{\ensuremath{\mathcal H}}

\newcommand{\addq}{\ensuremath{\oplus_q}}
\newcommand{\pol}{$\mathcal{F} = \{f_n(q)\}_{n=1}^{\infty}$}
\newcommand{\polg}{$\mathcal{G} = \{g_{n}(q)\}_{n=1}^{\infty}$}

\DeclareMathOperator{\Hom}{Hom}
\DeclareMathOperator{\End}{End}
\DeclareMathOperator{\supp}{supp}

\begin{document}

\vspace{3cm}

\title[Semidirect Products and Quantum Multiplication]{Semidirect Products and Functional Equations for Quantum Multiplication}
\author{Melvyn B. Nathanson}
\address{Department of Mathematics\\Lehman College (CUNY)\\Bronx, New York 10468}
\email{melvyn.nathanson@lehman.cuny.edu}
\thanks{This work was supported in part by grants from the NSA Mathematical Sciences Program
and the PSC-CUNY Research Award Program.}

\begin{abstract}
The quantum integer $[n]_q$ is the polynomial $1 + q + q^2 + \cdots + q^{n-1}$, and the sequence of polynomials $\{ [n]_q \}_{n=1}^{\infty}$ is a solution of the functional equation $f_{mn}(q) = f_m(q)f_n(q^m).$  In this paper, semidirect products of semigroups are used to produce families of functional equations that generalize the functional equation for quantum multiplication.
\end{abstract}

\date{\today}
\maketitle

\section{Multiplication of quantum integers}
Let \pol\ be a sequence of functions.  We define a binary operation $\otimes$
on the terms of the sequence \FF\ by 
\[
f_m(q)\otimes f_n(q) = f_m(q)f_n(q^m).
\]
For every positive integer $n$, the quantum integer $[n]_q$ is the polynomial
\[
[n]_q = 1+q+q^2+\cdots + q^{n-1}
\]
and
\[
[m]_q \otimes [n]_q = [mn]_q.
\]
Equivalently, the sequence of polynomials $\mathcal{F} = \{ [n]_q\}_{n=1}^{\infty}$ satisfies the functional equation
\beq        \label{semidirect:mfe}
 f_{mn}(q)  = f_m(q)f_n(q^m).
\eeq
It is an open problem to classify the sequences of functions (for example, polynomials, rational functions, or formal power series) 
that satisfy the functional equation~(\ref{semidirect:mfe}).
These problems have been studied in~\cite{bori-nath-wang04,nath03b,nath04c}, and related problems for addition of quantum integers in~\cite{kont-nath05,nath03d,nath06a}.
In this paper we apply the semidirect product of semigroups to produce families of functional equations that generalize the classical functional equation for quantum multiplication.

\section{Semidirect products of semigroups}
Let $S$ and $T$ be semigroups, written multiplicatively, 
with identity elements $e_S$ and $e_T$, respectively,
and let $\Hom(S,T)$ denote the set of all semigroup homomorphisms $\lambda:S\rightarrow T$
such that $\lambda(e_S) = e_T$.
In this paper, homomorphism always means semigroup homomorphism.
Let $\End(S) = \Hom(S,S)$ denote the semigroup of endomorphisms of $S$ 
under composition of maps.
Let 
\[
\alpha: T \rightarrow \End(S)
\]
be a homomorphism.  We denote the image of $t$ under $\alpha$ by $\alpha_t$.

We consider the set $S \times T$ with the
binary operation $\ltimes_{\alpha}$ defined as follows:
\[
(s_1,t_1) \ltimes_{\alpha} (s_2,t_2) = (s_1\alpha_{t_1}(s_2), t_1t_2).
\]
This multiplication is associative, since
\begin{align*}
(s_1,t_1)\ltimes_{\alpha} \left((s_2,t_2) \ltimes_{\alpha} (s_3,t_3) \right)
& = (s_1,t_1) \ltimes_{\alpha}(s_2\alpha_{t_2}(s_3), t_2t_3) \\
& = (s_1\alpha_{t_1}\left(s_2\alpha_{t_2}(s_3)\right), t_1(t_2t_3)) \\
& = (s_1\alpha_{t_1}(s_2)\alpha_{t_1}\alpha_{t_2}(s_3), (t_1t_2)t_3) \\
& = (s_1\alpha_{t_1}(s_2)\alpha_{t_1t_2}(s_3), (t_1t_2)t_3) \\
& = (s_1\alpha_{t_1}(s_2),t_1t_2)\ltimes_{\alpha} (s_3,t_3) \\
& = \left((s_1,t_1)\ltimes_{\alpha} (s_2,t_2)\right)\ltimes_{\alpha} (s_3,t_3).
\end{align*}
Since
\[
(e_S,e_T) \ltimes_{\alpha} (s,t) = (e_S\alpha_{e_T}(s),e_Tt) = (s,t)
\]
and
\[
(s,t)\ltimes_{\alpha}(e_S,e_T) = (s\alpha_t(e_S),te_T) = (s,t),
\]
it follows that the set $S \times T$ is a semigroup with identity $(e_S,e_T)$.
This semigroup is called 
the {\em semidirect product} of $S$ and $T$ with respect to $\alpha$, and denoted $
S \ltimes_{\alpha} T.$

\section{Multiplicative functional equations}
An {\em arithmetic function} is a complex-valued function whose domain is the multiplicative semigroup $\N = \{1,2,3,\ldots\}.$   An arithmetic function $u$ is called {\em completely multiplicative} if $u(mn) = u(m)u(n)$ for all $m,n \in \N$.  Then $\End(\N)$ is the semigroup of completely multiplicative arithmetic functions with identity $\epsilon$ defined by $\epsilon(n)=n$ for all $n \in \N.$

Let $A$ be an integral domain, that is, a commutative ring with identity and without zero divisors.  We apply the semidirect product construction in the case where $S = A[q]$  is the multiplicative semigroup of polynomials with coefficients in $A$, and $T = \N$ is the multiplicative semigroup of positive integers.  Let $\End^*(A[q])$ denote the subsemigroup of homomorphisms $\varphi: A[q]\rightarrow  A[q]$ such that $\varphi(0) = 0$ and $\varphi(f) \neq 0
$ for all $f \in A[q]$ with $f \neq 0$.  We fix a homomorphism 
\[
\alpha:\N \rightarrow \End^*(A[q])
\]
and consider the semidirect product $A[q]\ltimes_{\alpha} \N.$

For every homomorphism
\[
\Phi:\N \rightarrow A[q]\ltimes_{\alpha} \N,
\]
we define the sequence of polynomials
\pol\ and the arithmetic function $u: \N \rightarrow \N$ by
\[
\Phi(n) = (f_n(q),u(n)).
\]
We say that the sequence \FF\ and the function $u$ are \emph{associated} to $\Phi$.

For all positive integers $m$ and $n$ we have
\begin{align*}
(f_{mn}(q),u(mn)) & = \Phi(mn) \\ 
& = \Phi(m)\ltimes_{\alpha}\Phi(n) \\
& = (f_{m}(q),u(m))\ltimes_{\alpha}(f_{n}(q),u(n)) \\
& = \left(f_m(q)\alpha_{u(m)}(f_n(q)),u(m)u(n)\right).
\end{align*}
It follows that the sequence of polynomials \pol\ satisfies the functional equation
\beq        \label{semidirect:fe}
\boxed{f_{mn}(q) = f_m(q)\alpha_{u(m)}(f_n(q))}
\eeq
for all positive integers $m$ and $n$.
Moreover, $u(mn) = u(m)u(n)$ and so the arithmetic function $u$ is completely multiplicative.

Conversely, if $u \in \End(\N)$ and \pol\ is a sequence of polynomials 
that satisfies the functional equation~\eqref{semidirect:fe}, then the map
\[
\Phi:\N \rightarrow A[q]\ltimes_{\alpha} \N
\]
defined by
\[
\Phi(n) = (f_n(q),u(n))
\]
is a homomorphism, that is, $\Phi \in \Hom(\N,A[q]\ltimes_{\alpha} \N).$

For every completely multiplicative arithmetic function $u$ 
we denote by $\Hom(\alpha,u)$ the set of all homomorphisms
\[
\Phi:\N \rightarrow A[q]\ltimes_{\alpha} \N
\]
such that $u$ is associated to $\Phi$.
There is a one-to-one correspondence between polynomial solutions
of the functional equation~\eqref{semidirect:fe} and $\Hom(\alpha,u)$.  The general problem of quantum multiplication is to find all solutions of the functional equation~\eqref{semidirect:fe} for homomorphisms $\alpha$ and $u$.

\bt    \label{semidirect:theorem:abelian}
Let $\alpha:\N \rightarrow \End^*(A[q])$ and $u:\N \rightarrow \N$ be homomorphisms.  The set $\Hom(\alpha,u)$ is an abelian semigroup.
\et

\begin{proof}
If $\Phi_1, \Phi_2 \in \Hom(\alpha,u),$ then there exist sequences of polynomials
$\mathcal{F}_1 = \{f_{n,1}(q)\}_{n=1}^{\infty}$
and 
$\mathcal{F}_2 = \{f_{n,2}(q)\}_{n=1}^{\infty}$
such that
\[
\Phi_1(n) = (f_{n,1}(q),u(n)) \quad\text{ and }\quad \Phi_2(n) = (f_{n,2}(q),u(n))
\]
for all $n \in \N.$
Defining
\[
(\Phi_1\Phi_2)(n) = (f_{n,1}(q) f_{n,2}(q),u(n)),
\]
we obtain by a straightforward calculation that 
$\Phi_1 \Phi_2 \in \Hom(\alpha,u).$
This multiplication is associative with identity $n \mapsto (1,u(n))\in A[q]\times \N$ for all $n \in \N.$
\end{proof}

\bt
Let $\alpha:\N \rightarrow \End^*(A[q])$ and $u:\N \rightarrow \N$ be homomorphisms, and let
\[
\Phi:\N \rightarrow A[q]\ltimes_{\alpha} \N
\]
be a homomorphism in $\Hom(\alpha,u)$.  Let \pol\  be the sequence of polynomials associated to $\Phi.$   Then $\beta = \alpha \circ u: \N \rightarrow \End^*(A[q])$ is a homomorphism, and the map
\[
\Psi:\N\rightarrow A[q]\ltimes_{\beta} \N
\]
defined by 
\[
\Psi(n) = (f_n(q),n)
\]
is a homomorphism whose associated arithmetic function is the identity $\epsilon$.
Moreover,
\[
\Hom(\alpha,u) \cong \Hom(\alpha\circ u, \epsilon).
\]
\et

\begin{proof}
Let $\Phi \in \Hom(\alpha,u)$ and let \pol\  be the sequence of polynomials  associated to  $\Phi.$    Then  $\mathcal{F}$ satisfies the functional equation~(\ref{semidirect:fe}).  The composite function $\beta = \alpha\circ u$ is a homomorphism
\[
\beta:\N \rightarrow \End^*(A[q]).
\]
Define 
\[
\Psi:\N\rightarrow A[q]\ltimes_{\beta} \N
\]
by 
\[
\Psi(n) = (f_n(q),n).
\]
Then
\begin{align*}
\Psi(mn) & = (f_{mn}(q),mn) \\
& = (f_m(q)\alpha_{u(m)}(f_n(q)),mn)  \\
& = (f_m(q)\beta_{m}(f_n(q)),mn)  \\
& = (f_m(q)),m)\ltimes_{\beta}(f_n(q)),n)  \\
& = \Psi(m)\ltimes_{\beta}\Psi(n)
\end{align*}
and so $\Psi \in \Hom(\beta,\epsilon)$, that is, $\Psi$ is a homomorphism whose associated arithmetic function is the identity $\epsilon$.  Moreover, the map $\Phi \mapsto \Psi$ is a one-to-one homomorphism from $\Hom(\alpha,u)$ to $\Hom(\beta,\epsilon)$ with an inverse that maps $(f_n(q),n)$ to $(f_n(q),u(n))$.
\end{proof}

\section{Examples}
Let $u$ and $v$ be completely multiplicative arithmetic functions, and define
\[
\alpha: \N \rightarrow \End^*(A[q])
\]
by 
\[
\alpha_n(f)(q) = \left[f\left( q^{u(n)} \right)\right]^{v(n)}
\]
for all $n \in \N$.  We have $\alpha_n(f_1f_2)= \alpha_n(f_1)\alpha_n(f_2)$ for all $f_1,f_2 \in A[q]$, and $\alpha_n(f) = 0$ if and only if $f=0.$  Let $h(q) = \alpha_n(f)(q)$.  Then
\begin{align*}
\alpha_m \circ \alpha_n(f)(q) & = \alpha_m( \alpha_n(f) )(q)\\ 
& = \alpha_m(h)(q) \\
& = \left[ h\right(q^{u(m)} \left) \right]^{v(m)} \\
& = \left[ f \left( q^{u(n)u(m)} \right) \right]^{v(n)v(m)}\\
& = \left[ f\left( q^{u(mn)} \right) \right]^{v(mn)} \\
& = \alpha_{mn}(f)(q)
\end{align*}
and so $\alpha$ is a homomorphism.  In particular, for $v(n)=1$, the homomorphisms $\Phi:\N \rightarrow \End^*(A[q])\ltimes_{\alpha} \N$ with associated function $\epsilon$ are in one-to-one correspondence with the sequences \pol\ that satisfy the functional equation
\beq  \label{semidirect:ufe}
f_{mn}(q) = f_m(q)f_n\left(q^{u(m)}\right)
\eeq  
for all $m,n \in \N.$  This is called the \emph{twisted functional equation for quantum multiplication}.  There exist nontrivial solutions of this equation for every completely multiplicative function $u.$

\bt
Let $u:\N \rightarrow \N$ be a completely multiplicative.  The sequence
of quantum integers 
\[
\mathcal{F} = \{[u(n)]_q\}_{n=1}^{\infty}
\]
and the sequence  of powers
\[
\mathcal{F} = \{q^{u(n)-1}\}_{n=1}^{\infty}
\]
are solutions of the functional equation~\eqref{semidirect:ufe}.
\et

\begin{proof}
Let $f_n(q)=[u(n)]_q$.  Since the quantum integers satisfy the functional equation~\eqref{semidirect:mfe}, we have
\[
f_{mn}(q)  = [u(mn)]_q = [u(m)]_q [u(n)]_{q^{u(m)}} = f_m(q)f_n(q^{u(m)}).
\]
Similarly, if $f(n)=q^{u(n)-1}$, then 
\[
f_{mn}(q) = q^{u(mn)-1}=q^{u(m)-1}q^{u(m)(u(n)-1)} = f_m(q)f_n(q^{u(m)}).
\]
\end{proof}

\section{Support}
Let $\alpha:\N \rightarrow \End^*(A[q])$ and $\Phi:\N \rightarrow A[q]\ltimes_{\alpha} \N$ be homomorphisms and let \pol\ be the sequence of polynomials associated to $\Phi.$  We define the \emph{support} of  \FF\ and the \emph{support} of  $\Phi$ by
\[
\supp(\FF) = \supp(\Phi) = \{n \in \N : f_n(q) \neq 0\}.
\]
Since the sequence \pol\ satisfies the functional equation~\eqref{semidirect:fe}, it follows that if $m=n=1$, then $f_1(q)=f_1(q)^2$ and so  $f_1(q) = 0$ or $1$.  If $f_m(q) = 0,$ then $f_{mn}(q) = 0$ for all $n \in \N.$
In particular, if $f_1(q) = 0,$ then $f_n(q) = 0$ for all $n \in \N$ and $\supp(\FF) = \emptyset.$  
The sequence \FF\ is a called a \emph{nonzero solution} 
of the functional equation~(\ref{semidirect:fe}) if $\supp(\FF) \neq\emptyset.$  

Since $\alpha_m \in \End^*(A[q])$ for all $m \in \N$, it follows that $\alpha_m(f_n)(q) \neq 0$ if $f_n(q)\neq 0$.   The functional equation~\eqref{semidirect:fe} implies that if $f_m(q)$ and $f_n(q)$ are nonzero polynomials, then $f_{mn}(q)$ is a nonzero polynomial.  Conversely, if $f_{mn}(q) \neq 0$, then $f_m(q)\neq 0$ and $f_n(q)\neq 0.$  Thus, if \FF\ is a nonzero solution of the functional equation~(\ref{semidirect:fe}), then $\supp(\FF)$ is a subsemigroup of the positive integers.  Let $P$ be a set of prime numbers, and let $S(P)$ be the semigroup of positive integers all of whose prime factors belong to $P$.  If $P$ is the set of primes $p$ such that $f_p(q) \neq 0$, then $\supp(\FF) = S(P).$

For every homomorphism $\alpha:\N \rightarrow \End^*(A[q])$, for every completely multiplicative arithmetic function $u$,  
and for every set $P$ of prime numbers, 
we denote by $\Hom(\alpha,u,P)$ the set of all homomorphisms
\[
\Phi:\N \rightarrow A[q]\ltimes_{\alpha} \N
\]
such that $\Phi$ is associated to $u$ and $\supp(\Phi) = S(P)$.
The set $\Hom(\alpha,u,P)$ is also an abelian semigroup, and there is 
a one-to-one correspondence between solutions of the functional equation~\eqref{semidirect:fe} with support $S(P)$ and $\Hom(\alpha,u,P).$

\bt
For all homomorphisms $\alpha: \N \rightarrow \End^*(A[q])$ and $u: \N \rightarrow \N$ and for every set $P$ of prime numbers, the set $\Hom(\alpha,u,P)$ is a subsemigroup of  $\Hom(\alpha,u)$.
\et

\begin{proof}
This follows from Theorem~\ref{semidirect:theorem:abelian} and the fact that if $\FF_1$ and $\FF_2$ 
are solutions of the functional equation~\eqref{semidirect:fe} with $\supp(\FF_1) = \supp(\FF_2) = S(P)$, then the sequence $\FF_1\FF_2 = \{f_{n,1}(q)f_{n,2}(q)\}_{n=1}^{\infty}$
is also a solution of~(\ref{semidirect:fe}) with $\supp(\FF_1\FF_2) = \supp(\FF_1)\cap \supp(\FF_2) = S(P).$
\end{proof}

\bt
Let $P$ be a set of prime numbers, and let $\{h_p(q) : p \in P\}$
be a set of nonzero polynomials such that 
\beq      \label{semidirect:hfe}
h_{p_1}(q) \alpha_{u(p_1)}(h_{p_2}(q))
= h_{p_2}(q) \alpha_{u(p_2)}(h_{p_1}(q))
\eeq
for all $p_1,p_2 \in P.$  Then there exists a unique sequence \pol\
of polynomials that satisfies the functional equation~\eqref{semidirect:fe} such that 
$f_p(q) = h_p(q)$ for all $p \in P,$ and $\supp(\FF) = S(P).$
\et

\begin{proof}
Define $\varphi:P \rightarrow A[q]\ltimes_{\alpha} \N$ by $\varphi(p) = (h_p(q),u(p)).$
It follows from~(\ref{semidirect:hfe}) that 
$\varphi(p_1)\ltimes_{\alpha}\varphi(p_2) = \varphi(p_2)\ltimes_{\alpha}\varphi(p_1)$
for all $p_1,p_2 \in P,$ and so the map $\varphi$ extends to a unique homomorphism 
$\Phi$ from the free abelian semigroup $S(P)$ generated by $P$ into
$A[q]\ltimes_{\alpha} \N$.  Extend $\Phi$ to \N\ by setting $\Phi(n)=(0,u(n))$ for all $n \in \N \setminus S(P).$
This completes the proof.
\end{proof}

\section{Problems}
\benum
\item[(i)]
We already stated the general problem of quantum multiplication: Find all solutions of the functional equation~\eqref{semidirect:fe} for homomorphisms $\alpha$ and $u$.

\item[(ii)]
Let $u$ and $v$ be completely multiplicative arithmetic functions.
There is a homomorphism $\alpha:\N \rightarrow \End^*(A[q])$ defined by
\[
\alpha_n(f(q) = \left[f(q^{u(n)})\right]^{v(n)}
\]
for all $n \in \N.$  
Twisted quantum multiplication is the case $v(n)=1$, and ordinary quantum multiplication is the case $u(n) = n$ and $v(n)=1$.
Describe all homomorphisms $\alpha: \N \rightarrow \End^*(A[q])$.

\item[(iii)]
For every completely multiplicative function $u(n)$ and set $P$ of prime numbers,
classify the solutions of the twisted functional equation
\[
f_{mn}(q) = f_m(q)f_n(q^{u(m})
\] 
with support $S(P).$

\item[(iv)]
Borisov, Nathanson, and Wang~\cite{bori-nath-wang04} proved that the only solutions of the functional equation
\[
f_{mn}(q) = f_m(q)f_n(q^m)
\]
in rational functions with rational coefficients are essentially cyclotomic, that is, products of powers of quantum integers.
Are solutions of the twisted functional equation in rational functions with rational coefficients also cyclotomic?
\eenum

\bibliography{mathrefs,nathansn}
\end{document}